\documentclass[a4paper]{article}
\usepackage{amsfonts}
\usepackage{amsmath}

\input{tcilatex}
\begin{document}

\begin{eqnarray*}
&&\text{\textbf{The art of finding Calabi-Yau differential equations}} \\
&&\text{Dedicated to the 90-th birthday of Lars G\aa rding}
\end{eqnarray*}%
\begin{equation*}
\text{Gert Almkvist}
\end{equation*}

\textbf{1. Introduction.}

In 2003 van Straten and van Enckevort during a computer search found the
third order differential operator ($\eta )$ 
\begin{equation*}
L=\theta ^{3}-x(2\theta +1)(11\theta ^{2}+11\theta +5)+125x^{2}(\theta
+1)^{3}
\end{equation*}%
Here \ $\theta =x\frac{d}{dx}$. Then \ $Ly=0$ \ where 
\begin{equation*}
y=\sum_{n=0}^{\infty }A_{n}x^{n}
\end{equation*}%
for some unknown coefficients \ $A_{n},$ $n=0,1,2,...$It took me five years
before I found an explicit expression for $A_{n},$namely \ $A_{0}=1$ \ and%
\begin{equation*}
A_{n}=5\binom{2n}{n}^{-1}\binom{3n}{n}^{-1}\sum_{k=0}^{[n/5]}(-1)^{k}\frac{%
n-2k}{4n-5k}\binom{n}{k}^{-2}\binom{n}{5k}\binom{5n-5k}{n}^{-1}\frac{(5k)!}{%
k!^{5}}\frac{(5n-5k)!}{(n-k)!^{5}}
\end{equation*}%
for \ $n>0.$ How is it possible to find such a complicated formula? Let us
first point out that I immediately got an e-mail from Zudilin where he
(trivially) simplified the formula to%
\begin{equation*}
A_{n}=5\sum_{k=0}^{[n/5]}(-1)^{k}\frac{n-2k}{4n-5k}\binom{n}{k}^{3}\binom{%
4n-5k}{3n}
\end{equation*}%
and later (not so trivially) to 
\begin{equation*}
A_{n}=\sum_{k=0}^{[n/5]}(-1)^{k}\binom{n}{k}^{3}\left\{ \binom{4n-5k-1}{3n}+%
\binom{4n-5k}{3n}\right\} 
\end{equation*}%
Here it is evident that \ $A_{n}$ is an integer. But it was found in the
complicated form above and here we shall tell the long story how I came to
let Maple's "Zeilberger" find the recursion formula for this monster.

\textbf{2. Hypergeometric equations.}

The first Calabi-Yau differential equations connected to Calabi-Yau
manifolds needed in string theory appeared in the 1980-ies in the physics
literature. They were of the form%
\begin{equation*}
\theta ^{4}-xP(\theta )
\end{equation*}%
where \ $P(\theta )$ \ is a polynomial of degree \ $4.$ There were 13 of
them and I found the 14-th equation%
\begin{equation*}
\theta ^{4}-12^{6}x(\theta +\frac{1}{12})(\theta +\frac{5}{12})(\theta +%
\frac{7}{12})(\theta +\frac{11}{12})
\end{equation*}%
while writing the popular paper [1]. It was also found independently by
C.Doran and J.Morgan [14]. It is interesting that the same 14 4-tuples of
fractions ( like $(1/12,5/12,7/12,11/12)$) occur in J.Guillera's
Ramanujan-like formulas for \ $1/\pi ^{2}$ (see [14]). The hypergeometric
equations are \# 1-14 in [3] called the "Big Table" from now on.

\textbf{3. Formal definitions.}

After I had finished writing [1] I found the papers [10],[11] by Batyrev,
van Straten et al. There were another 14 equations (\# 15-28 in [3] ) coming
from geometry and they were not hypergeometric. What is the common property
of these 28 equations?

\textbf{Definition: }A Calabi-Yau differential equation is a 4-th order
differential equation with rational coefficients%
\begin{equation*}
y^{(4)}+a_{3}(x)y^{\prime \prime \prime }+a_{2}(x)y^{\prime \prime
}+a_{1}(x)y^{\prime }+a_{0}(x)y=0
\end{equation*}
satisfying the following conditions.

\textbf{1. }It is MUM (Maximal Unipotent Monodromy), i.e. the indicial
equation at $\ x=0$ has zero as a root of order 4. It means that there is a
Frobenius solution of the following form%
\begin{equation*}
y_{0}=1+A_{1}x+A_{2}x^{2}+...
\end{equation*}%
\begin{equation*}
y_{1}=y_{0}\log (x)+B_{1}x+B_{2}x^{2}+..
\end{equation*}%
\begin{equation*}
y_{2}=\frac{1}{2}y_{0}\log ^{2}(x)+(B_{1}x+B_{2}x^{2}+...)\log
(x)+C_{1}x+C_{2}x^{2}+...
\end{equation*}%
\begin{equation*}
y_{3}=\frac{1}{6}y_{0}\log ^{3}(x)+\frac{1}{2}(B_{1}x+B_{2}x^{2}+...)\log
^{2}(x)+(C_{1}x+C_{2}x^{2}+...)\log (x)+D_{1}x+D_{2}x^{2}+...
\end{equation*}%
It is very useful that Maple's "formal\_sol" produces the four solutions in
exactly this form (though labelled $\ 1-4$ )

\textbf{2. }The coefficents of the equation satisfy the identity%
\begin{equation*}
a_{1}=\frac{1}{2}a_{2}a_{3}-\frac{1}{8}a_{3}^{3}+a_{2}^{\prime }-\frac{3}{4}%
a_{3}a_{3}^{\prime }-\frac{1}{2}a_{3}^{\prime \prime }
\end{equation*}

\textbf{3.} Let \ $t=$ $y_{1}/y_{0}.$ Then%
\begin{equation*}
q=\exp (t)=x+c_{2}x^{2}+...
\end{equation*}%
can be solved%
\begin{equation*}
x=x(q)=q-c_{2}q^{2}+....
\end{equation*}%
which is called the "mirror map". We also construct the "Yukawa coupling"
defined by%
\begin{equation*}
K(q)=\frac{d^{2}}{dt^{2}}(\frac{y_{2}}{y_{0}})
\end{equation*}%
This can be expanded in a Lambert series%
\begin{equation*}
K(q)=1+\sum_{d=1}^{\infty }n_{d}\frac{d^{3}q^{d}}{1-q^{d}}
\end{equation*}%
where the \ $n_{d}$ \ are called "instanton numbers". For small \ $d$ \ the
\ $n_{d}$ \ are conjectured to count rational curves of degree \ $d$ \ on
the corresponding Calabi-Yau manifold. Then the third condition is

(\textbf{a}) $\ y_{0}$ has integer coefficients

(\textbf{b}) $q$ \ has integer coefficients

(\textbf{c}) There is a fixed integer \ $N_{0}$ \ such that all \ $%
N_{0}n_{d} $ \ are integers

\begin{equation*}
\end{equation*}

Actually it looks as if conditions \textbf{1} and \textbf{3a,b} would imply
conditions \textbf{2} and \textbf{3c}. At least we have not found any
counter example during six years of search.

\begin{equation*}
\end{equation*}

\textbf{3. Pullbacks of 5-th order equations.}

The condition \textbf{2} \ is equivalent to%
\begin{equation*}
\begin{vmatrix}
y_{0} & y_{3} \\ 
y_{0}^{\prime } & y_{3}^{\prime }%
\end{vmatrix}%
=%
\begin{vmatrix}
y_{1} & y_{2} \\ 
y_{1}^{\prime } & y_{2}^{\prime }%
\end{vmatrix}%
\end{equation*}%
This means that the six wronskians formed by the four solutions to our
Calabi-Yau equation reduce to five. Hence they satisfy a 5-th order
differential equation%
\begin{equation*}
w^{(5)}+b_{4}w^{(4)}+b_{3}w^{\prime \prime \prime }+b_{2}w^{\prime \prime
}+b_{1}w^{\prime }+b_{0}w=0
\end{equation*}%
The condition 2 for the 4-th order equation leads to a corresponding
condition for the 5-th order equation

\textbf{2}$_{\text{\textbf{5}}}$%
\begin{equation*}
b_{2}=\frac{3}{5}b_{3}b_{4}-\frac{4}{25}b_{4}^{3}+\frac{3}{2}b_{3}^{\prime }-%
\frac{6}{5}b_{4}b_{4}^{\prime }-b_{4}^{\prime \prime }
\end{equation*}

The special equation (coming from number theory, see Zudilin [22] )%
\begin{equation*}
\theta ^{5}-3x(2\theta +1)(3\theta ^{2}+3\theta +1)(15\theta ^{2}+15\theta
+4)-3x^{2}(\theta +1)^{3}(3\theta +2)(3\theta +4)
\end{equation*}%
satisfies this equation and it is possible to find its 4-th order "pullback"
(\# 32 in the Big Table). It has degree \ $8$ since%
\begin{equation*}
\theta ^{4}+...+3^{8}x^{8}(3\theta +1)^{2}(3\theta +2)^{2}
\end{equation*}%
with large coefficients. We say that two Calabi-Yau equations are equivalent
if they have the same instanton numbers. This is the case under the
transformation%
\begin{equation*}
Y(x)=f(x)y(g(x))
\end{equation*}%
where%
\begin{equation*}
g(x)=x+a_{2}x^{2}+...
\end{equation*}%
Conversely it is conjectured that equivalence implies the existence of
algebraic \ $f(x)$ \ and \ $g(x)$ \ as above (see [4]).Yifan Yang ([21])
suggested a different but equivalent pullback which usually cuts the degree
in half. It also has a symmetry reducing the number of coefficient another
50\% (see [6],[7])..It depends on the following. Let%
\begin{equation*}
w_{0}=x%
\begin{vmatrix}
y_{0} & y_{1} \\ 
y_{0}^{\prime } & y_{1}^{\prime }%
\end{vmatrix}%
\end{equation*}%
\begin{equation*}
w_{1}=x%
\begin{vmatrix}
y_{0} & y_{2} \\ 
y_{0}^{\prime } & y_{2}^{\prime }%
\end{vmatrix}%
\end{equation*}%
Then we have the following identity, "The double wronskian is almost the
square"%
\begin{equation*}
\begin{vmatrix}
w_{0} & w_{1} \\ 
w_{0}^{\prime } & w_{1}^{\prime }%
\end{vmatrix}%
=x^{2}y_{0}^{2}\exp (-\frac{1}{2}\int a_{3}dx)
\end{equation*}%
(for a proof see [6]). Thus \ $y_{0}$ \ is, up to a factor, equal to the
square root of the wronskian of \ $w_{0}$ \ and \ $w_{1}.$ We give the
Yifan-Yang pullback of \# 32%
\begin{equation*}
\theta ^{4}-x\left\{ 540(\theta +\frac{1}{2})^{4}+486(\theta +\frac{1}{2}%
)^{2}+\frac{57}{4}\right\}
\end{equation*}%
\begin{equation*}
+x^{2}\left\{ 72846(\theta +1)^{4}+\frac{6915}{2}(\theta +1)^{2}+\frac{3}{4}%
\right\}
\end{equation*}%
\begin{equation*}
+x^{3}\left\{ 14580(\theta +\frac{3}{2})^{4}+12717(\theta +\frac{3}{2}%
)^{2}+324\right\}
\end{equation*}%
\begin{equation*}
+\frac{9}{16}x^{4}(6\theta +11)^{2}(6\theta +13)^{2}
\end{equation*}%
The most striking result from the Yifan-Ying pullback are 14 degree 2
equations coming from hypergeometric 5-th order equations. We show the
following example. Consider%
\begin{equation*}
\theta ^{5}-4\cdot 12^{6}x(\theta +\frac{1}{2})(\theta +\frac{1}{12})(\theta
+\frac{5}{12})(\theta +\frac{7}{12})(\theta +\frac{11}{12})
\end{equation*}%
which has the Y-Y pullback%
\begin{equation*}
\theta ^{4}-144x(165888\theta ^{4}+331776\theta ^{3}+386496\theta
^{2}+220608\theta +47711)
\end{equation*}%
\begin{equation*}
+2^{22}3^{10}x^{2}(4\theta +3)(4\theta +4)(6\theta +5)(6\theta +7)
\end{equation*}%
an equation not found in the computer search mentioned in the introduction
because the numbers are too big.

But unfortunately 5-th order differential equations satisfying \ \textbf{2}$%
_{\text{\textbf{5 \ \ }}}$are rare (except the ones constructed from known
4-th order C-Y equations). In [7] they are listed. In the new edition of the
Big Table there are two new ones, \# 355, 356.%
\begin{equation*}
\end{equation*}

\textbf{4. Using Maple for finding and factoring differential equations.}

Using Maple's "Zeilberger" is the best way to find Calabi-Yau differential
equations coming from simple sums of \ products of binomial coefficients.

\textbf{Example: \#15. }Consider the sum%
\begin{equation*}
A_{n}=\frac{(3n)!}{n!^{3}}\sum_{k=0}^{n}\binom{n}{k}^{3}
\end{equation*}%
Then Maple finds the recursion formula for \ $A_{n}$. Here \ $NA_{n}=A_{n+1}$%
\ \ 

with(SumTools[Hypergeometric]);

r:=Zeilberger((3n)!/n!\symbol{94}3*binomial(n,k)\symbol{94}3,n,k,N)[1];

$r:=(n+2)^{4}N^{2}-3\left\{ 3(n+1)+1)\right\} \left\{ 3(n+1)+2\right\}
\left\{ 7(n+1)^{2}+7(n+1)+2\right\} N$

$-72(3n+1)(3n+2)(3n+4)(3n+5);$

Maple finds the recursion of lowest order (degree in \ $N$ ), say%
\begin{equation*}
c_{0}(n)N^{p}+c_{1}(n)N^{p-1}+...+c_{p}(n)
\end{equation*}%
Converting to a differential operator we find%
\begin{equation*}
c_{0}(\theta -p)+xc_{1}(\theta -p+1)+...+x^{p}c_{p}(\theta )
\end{equation*}%
Hence the C-Y condition \textbf{1 }(MUM) is equivalent to%
\begin{equation*}
c_{0}(n)=(n+p)^{4}
\end{equation*}%
But this is not always the case.

\textbf{Example: \#22.} Consider the sum%
\begin{equation*}
A_{n}=\sum_{k=0}^{n}\binom{n}{k}^{5}
\end{equation*}%
Then we compute

r:=Zeilberger(binomial(n,k)\symbol{94}5,n,k,N)[1];

$r:=c_{0}(n)N^{3}+c_{1}(n)N^{2}+c_{2}(n)N+c_{3}(n);$

where%
\begin{equation*}
c_{0}(n)=(55n^{2}+143n+94)(n+3)^{4}
\end{equation*}%
This will give a differential equation \ $L$ \ of degree 6 which can be
factored in Maple. Observe that first we have to convert the differential
equation so that \ $\theta $ \ is replaced by \ $\frac{d}{dx}.$ This can be
done in Maple using Stirling numbers (see p.19 in [6]).

with(DEtools);

DFactor(L,[Dx,x], `one step `);

L2:=\%[2];

Here L2 is the right factor of L.

This strategy usually does not work if \ $c_{0}(n)/(n+p)^{4}$ \ contains
only linear factors of \ the form \ $n+a,2n+a,3n+a,4n+a.$ There are two
exceptions:

\#354 \ $c_{0}(n)=(n+1)(7n+12)(n+3)^{4}$

\#361 \ $c_{0}(n)=$\ $(2n+5)(8n+77)(n+4)^{4}$

We have also found three cases where we have an irreducible factor of degree
four.

\#251 \ $c_{0}(n)=(441n^{4}+3780n^{3}+11634n^{2}+15207n+7202)(n+4)^{4}$

\#299 \ $c_{0}(n)=$ $%
(4704n^{4}+40572n^{3}+117732n^{2}+133251n+49684)(n+4)^{4}$

\#367 $\ c_{0}(n)=(1888n^{4}-13592n^{3}+35556n^{2}-39931n+16322)(n+4)^{4}$

There is also a case with an reducible factor of degree 3

\#235 \ $c_{0}(n)=(5n+9)(6n^{2}+22n+19)(n+4)^{4}$\ 
\begin{equation*}
\end{equation*}

\textbf{5. Multiple sums of binomial coefficients.}

There is a "MultiZeilberger" but it is too slow to be of any practical use.
So we use "brute force" instead. Let us take an

\textbf{Example. \#349.} Let%
\begin{equation*}
A_{n}=\sum_{k,l}(-1)^{n+k}3^{n-3k}\binom{n}{3k}\frac{(3k)!}{k!^{3}}\binom{n}{%
l}\binom{2k}{n-l}\binom{2l}{n}
\end{equation*}%
Assume that the annihilating differential equation is%
\begin{equation*}
Ly=q_{0}y^{(4)}+q_{1}y^{\prime \prime \prime }+q_{2}y^{\prime \prime
}+q_{1}y^{\prime }+q_{0}y=0
\end{equation*}%
where%
\begin{equation*}
q_{0}=a_{4}x^{4}+a_{5}x^{5}+...+a_{16}x^{16}
\end{equation*}%
\begin{equation*}
q_{1}=b_{3}x^{3}+b_{4}x^{4}+...+b_{15}x^{15}
\end{equation*}%
\begin{equation*}
q_{2}=c_{2}x^{2}+c_{3}x^{3}+...+c_{14}x^{14}
\end{equation*}%
\begin{equation*}
q_{3}=d_{1}x+d_{2}x^{2}+...+d_{13}x^{13}
\end{equation*}%
\begin{equation*}
q_{4}=e_{0}+e_{1}x+...+e_{12}x^{12}
\end{equation*}%
with 65 unknown coefficients \ $a_{4},..e_{12}.$ Then we compute 70
coefficients of\ $A_{n}$ (takes about one second)\ and substituting%
\begin{equation*}
y=\sum_{n=0}^{69}A_{n}x^{n}
\end{equation*}%
in \ $Ly$ \ we get 70 linear equations which is solved in a few seconds on a
laptop. \#349 is the only known case where we really need coefficients of
degree 12.%
\begin{equation*}
\end{equation*}

\textbf{6. Hadamard products.}

If%
\begin{equation*}
u=\sum_{n=0}^{\infty }b_{n}x^{n}\text{ \ \ }v=\sum_{n=0}^{\infty }c_{n}x^{n}
\end{equation*}%
are two D-finite (i.e. satisfying differential equations with polynomial
coefficients) power series, then the Hadamard product 
\begin{equation*}
y=u\ast v=\sum_{n=0}^{\infty }a_{n}x^{n}=\sum_{n=0}^{\infty }b_{n}c_{n}x^{n}
\end{equation*}%
is also D-finite. It was suggested by Duco van Straten that if \ $u$ \ and \ 
$v$ \ satisfied "nice" second order equations then \ $u\ast v$ \ could
satisfy a Calabi-Yau equation. This was indeed the case for the following
type%
\begin{equation*}
\theta ^{2}-x(a\theta ^{2}+a\theta +b)+cx^{2}(\theta +1)^{2}
\end{equation*}%
with integer \ $a,b,c$, a class of differential equations studied by Don
Zagier ([21]). There are 10 such equations listed \ \ as \ (a),(b),..,(j) \
in [2] . There are also 10 third order equations 
\begin{equation*}
\theta ^{3}-x(2\theta +1)(\widehat{a}\theta ^{2}+\widehat{a}\theta +\widehat{%
b})+\widehat{c}x^{2}(\theta +1)^{3}
\end{equation*}%
suitable for Hadamard products with e.g. $\binom{2n}{n}^{2}$ \ giving 5-th
order equations. There are relations between the coefficients \ $a,b,c$ \
and \ $\widehat{a},\widehat{b},\widehat{c}$ \ described in [4] \ There are
many equivalences between the various Hadamard products. All this treated in
[4].%
\begin{equation*}
\end{equation*}

\textbf{7. The mirror at infinity.}

There are many Calabi-Yau equations ending with a term \ $cx^{p}(\theta
+1)^{4}$. For all these equations there is a "mirror at infinity" obtained
by the transformation \ $\theta \longrightarrow -\theta -1$ \ and \ $%
x\longrightarrow ax^{-1}$ for a suitable constant \ $a.$ E. R\"{o}dland did
this for \# 27 in [19].

\textbf{Example. \# 193. }We have%
\begin{equation*}
A_{n}=\sum_{k,l}\binom{n}{k}^{2}\binom{n}{l}^{2}\binom{k+l}{l}\binom{n+k+l}{n%
}
\end{equation*}%
and the equation%
\begin{equation*}
7^{2}\theta ^{4}-7x(1135\theta ^{4}+2204\theta ^{3}+1683\theta
^{2}+581\theta +77)
\end{equation*}%
\begin{equation*}
+x^{2}(28723\theta ^{4}+40708\theta ^{3}+13260\theta ^{2}-1337\theta -896)
\end{equation*}%
\begin{equation*}
-x^{3}(32126\theta ^{4}+38514\theta ^{3}+26511\theta ^{2}+10731\theta +1806)
\end{equation*}%
\begin{equation*}
+7\cdot 11x^{4}(130\theta ^{4}+254\theta ^{3}+192\theta ^{2}+65\theta
+8)+11^{2}x^{5}(\theta +1)^{4}
\end{equation*}%
We make the substitutions \ $\theta \longrightarrow -\theta -1$ and \ $%
x\longrightarrow 3^{-4}x^{-1}$. It follows%
\begin{equation*}
11^{2}\theta ^{4}-7\cdot 11x(130\theta ^{4}+266\theta ^{3}+210\theta
^{2}+77\theta +11)
\end{equation*}%
\begin{equation*}
-x^{2}(32126\theta ^{4}+89990\theta ^{3}+103725\theta ^{2}+55253\theta
+11198)
\end{equation*}%
\begin{equation*}
+x^{3}(28723\theta ^{4}+74184\theta ^{3}+63474\theta ^{2}+20625\theta +1716)
\end{equation*}%
\begin{equation*}
-7x^{4}(1135\theta ^{4}+2336\theta ^{3}+1881\theta ^{2}+713\theta
+110)+7^{2}x^{5}(\theta +1)^{4}
\end{equation*}%
which is \#198. In this case we also know a formula for the coefficients%
\begin{equation*}
A_{n}=\sum_{k,l}\binom{n}{k}^{2}\binom{n}{l}^{2}\binom{k+l}{l}\binom{2n-k}{n}
\end{equation*}%
which is rather unusual (11 out of 36 cases). Also equations ending with \ $%
cx^{p}(2\theta +1)^{4}$ can be treated in a similar way by the substitutions
\ $\theta \longrightarrow -\theta -1/2$ \ and \ $x\longrightarrow ax^{-1}.$%
\begin{equation*}
\end{equation*}

\textbf{8.Harmonic Sums.}

At the end of the paper [18] by P.Paule and C.Schneider there is a remark
that using "Zeilberger" on%
\begin{equation*}
"A_{n}"=\sum_{k}(n-2k)\binom{n}{k}^{7}
\end{equation*}%
which is identically zero by symmetry, one obtains the same recursion
formula as they obtained for%
\begin{equation*}
A_{n}=\sum_{k}\binom{n}{k}^{7}\left\{ 1+k(-7H_{k}+7H_{n-k})\right\}
\end{equation*}%
where 
\begin{equation*}
H_{n}=\sum_{j=1}^{n}\frac{1}{j}
\end{equation*}%
if \ $n\geq 1$ \ and \ $H_{n}=0$ \ if \ $n\leq 0.$ The recursion gives the
differential equation \# 27 in the Big Table.

\textbf{Lemma 1.} We have 
\begin{equation*}
\frac{d}{dn}n!=n!(H_{n}-\gamma )
\end{equation*}%
where \ $\gamma $ \ is Euler's constant.

Using this we find%
\begin{equation*}
-\frac{1}{2}\frac{d}{dk}"A_{n}"=A_{n}
\end{equation*}%
Indeed%
\begin{equation*}
-\frac{1}{2}\frac{d}{dk}\sum_{k=0}^{n}(n-2k)\binom{n}{k}^{7}=\frac{d}{dk}%
\sum_{k=0}^{n}k\binom{n}{k}^{7}=\sum_{k=0}^{n}\binom{n}{k}^{7}\left\{
1+k(-7H_{k}+7H_{n-k})\right\}
\end{equation*}

In this way 28 equations of type 
\begin{equation*}
"A_{n}"=\sum_{k}(n-2k)C(n,k)
\end{equation*}%
with \ 
\begin{equation*}
C(n,n-k)=C(n,k)
\end{equation*}%
were found, the last being \#360 with%
\begin{equation*}
A_{n}=\sum_{k}\binom{n}{k}\binom{n+3k}{n}\binom{4n-3k}{n}\frac{(3k)!}{k!^{3}}%
\frac{(3n-3k)!}{(n-k)!^{3}}
\end{equation*}%
\begin{equation*}
\left\{ 1+k(-4H_{k}+4H_{n-k}+3H_{n+3k}-3H_{4n-3k})\right\}
\end{equation*}

But sometimes it is not enough to take the derivative of \ $"A_{n}"$, we
also have to sum over negative \ $k.$(this was pointed out to me by
Christian Krattenthaler). For this we need

\textbf{Lemma 2. }Let \ $n$ \ be a positive integer. Then 
\begin{equation*}
\Gamma (-n+x)=\frac{(-1)^{n}}{n!}x^{-1}+O(1)
\end{equation*}%
when \ $x\rightarrow 0$.

\textbf{Proof. }We have%
\begin{equation*}
\Gamma (t)\Gamma (1-t)=\frac{\pi }{\sin (\pi t)}
\end{equation*}%
If \ $t=n+1-x$ \ we obtain%
\begin{equation*}
\Gamma (-n+x)=\frac{\pi }{\sin (\pi (n+1-x))}\frac{1}{\Gamma (n+1-x)}
\end{equation*}%
\begin{equation*}
=-\frac{\pi (-1)^{n+1}}{\sin (\pi x)}\frac{1}{\Gamma (n+1-x)}=\frac{(-1)^{n}%
}{x\Gamma (n+1)}+O(1)
\end{equation*}%
To illustrate this we consider

\textbf{Example \#264. }Let%
\begin{equation*}
"A_{n}"=16^{-n}\binom{2n}{n}^{2}\sum_{k}(n-2k)\binom{n}{k}\binom{2k}{k}%
\binom{2n-2k}{n-k}\binom{2n+2k}{n+k}^{2}\binom{4n-2k}{2n-k}^{2}\binom{2n}{k}%
^{-1}\binom{2n}{n-k}^{-1}
\end{equation*}%
Using \textbf{Lemma 2} we compute%
\begin{equation*}
\binom{n}{-k-\varepsilon }=\frac{(-1)^{k}}{k}\binom{n+k}{n}^{-1}\varepsilon
+O(\varepsilon ^{2})
\end{equation*}%
\begin{equation*}
\binom{-2k-2\varepsilon }{-k-\varepsilon }=\frac{1}{k}\binom{2k}{k}%
^{-1}\varepsilon +O(\varepsilon ^{2})
\end{equation*}%
\begin{equation*}
\binom{2n}{-k-\varepsilon }=\frac{(-1)^{k}}{k}\binom{2n+k}{2n}%
^{-1}\varepsilon +O(\varepsilon ^{2})
\end{equation*}%
Collecting this we find \ the derivative of \ $"A_{n}"$ at \ $-k$ \ and the
correct formula 
\begin{equation*}
A_{n}=16^{-n}\binom{2n}{n}^{2}\left\{ 
\begin{array}{c}
\sum_{k=0}^{n}\binom{n}{k}\binom{2k}{k}\binom{2n-2k}{n-k}\binom{2n+2k}{n+k}%
^{2}\binom{4n-2k}{2n-k}^{2}\binom{2n}{k}^{-1}\binom{2n}{n-k}^{-1} \\ 
\times \left\{
1+k(-2H_{k}+2H_{n-k}-3H_{n+k}+3H_{2n-k}+2H_{2k}-2H_{2n-2k}+4H_{2n+2k}-4H_{4n-2k})\right\} 
\\ 
+\dsum\limits_{k=1}^{n}\frac{n+2k}{k}\binom{2n+k}{2n}\binom{2n+2k}{n+k}%
\binom{2n-2k}{n-k}^{2}\binom{4n+2k}{2n+k}^{2}\binom{2k}{k}^{-1}\binom{n+k}{n}%
^{-1}\binom{2n}{n+k}^{-1}%
\end{array}%
\right\} 
\end{equation*}

\bigskip

We have a family of equations with%
\begin{equation*}
"A_{n}"=2^{-sn}\sum_{k}(n-2k)\binom{n}{k}^{7-2p}\binom{2k}{k}^{p}\binom{2n-2k%
}{n-k}^{p}
\end{equation*}%
with%
\begin{equation*}
\begin{tabular}{|l|l|l|l|}
\hline
p & s & \# & \#$^{\infty }$ \\ \hline
0 & 0 & 27 & 243 \\ \hline
1 & 0 & 212 & 117 \\ \hline
2 & 0 & 246 & 247 \\ \hline
3 & 0 & $\sim $6* &  \\ \hline
4 & 0 & $\sim $6* &  \\ \hline
5 & 2 & 247 & 246 \\ \hline
6 & 8 & 117 & 212 \\ \hline
7 & 14 & 243 & 27 \\ \hline
\end{tabular}%
\end{equation*}

\textbf{.}

For \ \#117 \ the technique used above \ gives%
\begin{equation*}
A_{n}=256^{-n}\left\{ 
\begin{array}{c}
\sum_{k=0}^{n}\binom{n}{k}^{-5}\binom{2k}{k}^{6}\binom{2n-2k}{n-k}%
^{6}\left\{ 1+k(-7H_{k}+7H_{n-k}+12H_{2k}-12H_{2n-2k})\right\} \\ 
+\sum_{k=1}^{\infty }(-1)^{k}\frac{n+2k}{k}\binom{n+k}{n}^{5}\binom{2k}{k}%
^{-6}\binom{2n+2k}{n+k}^{6}%
\end{array}%
\right\}
\end{equation*}%
which is nonsense since we cannot sum the infinite sum (maybe PARI can do
it?). Indeed let \ $n=1$ \ in the second sum. Then 
\begin{equation*}
256^{-1}\sum_{k=1}^{\infty }(-1)^{k}\frac{1+2k}{k}\binom{1+k}{k}^{5}\binom{2k%
}{k}^{-6}\binom{2+2k}{1+k}^{6}x^{k})=-\frac{13}{4}+(-\frac{57}{4}-\frac{1}{4}%
\log (2))(x-1)+O((x-1)^{2})
\end{equation*}%
where the constant term exactly cancels the value of the first sum. The same
occurs when \ $n=2$ \ and \ $n=3$. So there is still no formula known for \ $%
A_{n}$ for \# 117. There are 12 other cases (like \#243) where taking the
derivative with respect to \ $k$ \ does not work, but for which there exist
other formulas for \ $A_{n}.$

\bigskip

\textbf{9, Empty sums.}

\textbf{Example \# 133. }Consider the sum%
\begin{equation*}
"A_{n}"=\binom{2n}{n}^{2}\sum_{k}(n-2k)\binom{n}{k}^{-1}\binom{n}{3k}\binom{n%
}{3n-3k}\frac{(3k)!}{k!^{3}}\frac{(3n-3k)!}{(n-k)!^{3}}
\end{equation*}%
"Zeilberger" gives the recursion%
\begin{equation*}
(n+2)^{4}N^{2}-12(2n+3)^{2}(3n^{2}+9n+7)N+432(2n+1)^{2}(2n+3)^{2}
\end{equation*}%
which corresponds to the differential equation%
\begin{equation*}
\theta ^{4}-12x(2\theta +1)^{2}(3\theta ^{2}+3\theta +1)+432x^{2}(2\theta
+1)^{2}(2\theta +3)^{2}
\end{equation*}%
which we recognize as the Hadamard product \ \ $A\ast f$ (see [2]).

Looking at \ $"A_{n}"$ \ we observe that in order to \ $\binom{n}{3k}$ to be
nonzero we need \ $k\leq n/3.$ Similarly \ $\binom{n}{3n-3k}$ \ is nonzero
only for \ $k\geq 2n/3.$ So the sum is not only zero, it is also empty.
Consider the case \ $k\leq n/3$. Then%
\begin{equation*}
\binom{n}{3n-3k-3\varepsilon }=\frac{n!}{\Gamma (3n-3k-3\varepsilon
+1)\Gamma (-(2n-3k-1)-3\varepsilon )}
\end{equation*}%
\begin{equation*}
=\frac{n!}{(3n-3k)!}(\frac{(-1)^{2n-3k-1}}{(2n-3k-1)!}\frac{1}{%
(-3\varepsilon )})^{-1}=\frac{3(-1)^{k}}{2n-3k}\binom{3n-3k}{n}%
^{-1}\varepsilon +O(\varepsilon ^{2})
\end{equation*}%
which gives 
\begin{equation*}
A_{n}=3\binom{2n}{n}^{2}\sum_{k=0}^{[n/3]}(-1)^{k}\frac{n-2k}{2n-3k}\binom{n%
}{k}^{-1}\binom{n}{3k}\binom{3n-3k}{n}^{-1}\frac{(3k)!}{k!^{3}}\frac{(3n-3k)!%
}{(n-k)!^{3}}
\end{equation*}%
which simplifies to%
\begin{equation*}
A_{n}=3\binom{2n}{n}^{2}\sum_{k=0}^{[n/3]}(-1)^{k}\frac{n-2k}{2n-3k}\binom{n%
}{k}^{2}\binom{2n-3k}{n}
\end{equation*}%
for n\TEXTsymbol{>}0.

This is just a special case of%
\begin{equation*}
"A_{n}"=C_{n}\sum_{k}(n-2k)\binom{n}{k}^{a}E(n,k)\binom{n}{3k}\binom{n}{3n-3k%
}\frac{(3k)!}{k!^{3}}\frac{(3n-3k)!}{(n-k)!^{3}}
\end{equation*}%
with the following table%
\begin{equation*}
\begin{tabular}{|l|l|l|l|}
\hline
\# & $a$ & $C_{n}$ & $E(n,k)$ \\ \hline
133 & $-1$ & $\binom{2n}{n}^{2}$ & $1$ \\ \hline
279 & $1$ & $1$ & $1$ \\ \hline
334 & $1$ & $\binom{2n}{n}^{2}$ & $\binom{2n}{k}^{-1}\binom{2n}{n-k}^{-1}$
\\ \hline
281 & $2$ & $\binom{2n}{n}^{2}$ & $\binom{2n}{k}^{-1}\binom{2n}{n-k}^{-1}$
\\ \hline
363 & $-3$ & $\binom{2n}{n}\binom{4n}{2n}$ & $\binom{2n}{k}\binom{2n}{n-k}$
\\ \hline
352 & $-1$ & $1$ & $\binom{2k}{k}\binom{2n-2k}{n-k}$ \\ \hline
253 & $0$ & $\binom{2n}{n}$ & $\binom{2k}{k}\binom{2n-2k}{n-k}\binom{2n}{2k}%
^{-1}$ \\ \hline
353 & $-1$ & $\binom{2n}{n}$ & $\binom{2k}{k}\binom{2n-2k}{n-k}\binom{2n}{2k}%
^{-1}$ \\ \hline
350 & $-2$ & $\binom{2n}{n}^{3}$ & $\binom{2k}{k}\binom{2n-2k}{n-k}\binom{n+k%
}{n}^{-1}\binom{2n-k}{n}^{-1}$ \\ \hline
\end{tabular}%
\end{equation*}%
The computations before show that 
\begin{equation*}
A_{n}=3C_{n}\sum_{k=0}^{[n/3]}(-1)^{k}\frac{n-2k}{2n-3k}\binom{n}{k}%
^{3+a}E(n,k)\binom{2n-3k}{n}
\end{equation*}%
for \ $n\geq 1.$

Replacing \ $3$ \ by \ $4$ \ we get%
\begin{equation*}
"A_{n}"=C_{n}\sum_{k}(n-2k)\binom{n}{k}^{a}E(n,k)\binom{n}{4k}\binom{n}{4n-4k%
}\frac{(4k)!}{k!^{4}}\frac{(4n-4k)!}{(n-k)!^{4}}
\end{equation*}%
with the table%
\begin{equation*}
\begin{tabular}{|l|l|l|l|}
\hline
\# & $a$ & $C_{n}$ & $E(n,k)$ \\ \hline
300 & $-2$ & $\frac{(5n)!}{(2n)!n!^{3}}$ & $\binom{n+k}{n}^{-1}\binom{2n-k}{n%
}^{-1}$ \\ \hline
36 & $-2$ & $\binom{2n}{n}$ & $1$ \\ \hline
364 & $-1$ & $1$ & $1$ \\ \hline
357 & $0$ & $\binom{2n}{n}^{-1}$ & $1$ \\ \hline
205 & $0$ & $1$ & $\binom{2n}{2k}^{-1}$ \\ \hline
365 & $-3$ & $1$ & $\binom{2k}{k}\binom{2n-2k}{n-k}$ \\ \hline
\end{tabular}%
\end{equation*}%
Using \textbf{Lemma 2} one easily shows%
\begin{equation*}
A_{n}=4\binom{2n}{n}C_{n}\sum_{k=0}^{[n/4]}\frac{n-2k}{3n-4k}\binom{n}{k}%
^{4+a}E(n,k)\binom{3n-4k}{2n}
\end{equation*}

Replacing \ $4$ \ by \ $5$ \ we get%
\begin{equation*}
"A_{n}"=C_{n}\sum_{k}(n-2k)\binom{n}{k}^{a}E(n,k)\binom{n}{5k}\binom{n}{5n-5k%
}\frac{(5k)!}{k!^{5}}\frac{(5n-5k)!}{(n-k)!^{5}}
\end{equation*}%
with the table%
\begin{equation*}
\begin{tabular}{|l|l|l|l|}
\hline
\# & $a$ & $C_{n}$ & $E(n,k)$ \\ \hline
354 & $-3$ & $1$ & $1$ \\ \hline
$B\ast \eta $ & $-2$ & $1$ & $1$ \\ \hline
\end{tabular}%
\end{equation*}%
We get%
\begin{equation*}
A_{n}=5\binom{2n}{n}\binom{3n}{n}C_{n}\sum_{k=0}^{[n/5]}(-1)^{k}\frac{n-2k}{%
4n-5k}\binom{n}{k}^{5+a}E(n,k)\binom{4n-5k}{3n}
\end{equation*}%
for \ $n>0.$ For the case \ $B\ast \eta $ \ we get since \ $B\ast $ \ is
multiplication by \ $\frac{(3n)!}{n!^{3}}=\binom{2n}{n}\binom{3n}{n}$ the
formula for the coefficient for the third order equation \ $\eta $%
\begin{equation*}
A_{n}=5\sum_{k=0}^{[n/5]}(-1)^{k}\frac{n-2k}{4n-5k}\binom{n}{k}^{2}\binom{%
4n-5k}{3n}
\end{equation*}%
the equation mentioned in the Introduction.

Finally there is the case \#347 with%
\begin{equation*}
A_{n}=6\binom{2n}{n}^{2}\sum_{k=0}^{[n/6]}\frac{n-2k}{5n-6k}\binom{n}{k}^{2}%
\binom{5n-6k}{3n}
\end{equation*}

One of the most intricate cases is \#305 with%
\begin{equation*}
"A_{n}"=\binom{2n}{n}^{2}\sum_{k}(n-2k)\binom{n+2k}{k}\binom{3n-2k}{n-k}%
\binom{2n+4k}{n+2k}\binom{6n-4k}{3n-2k}\binom{3n}{n+k}
\end{equation*}%
There one has to have different sums for \ $-k$ \ depending on if \ $k<n/2$
\ or \ $k>n/2.$ See the final result in its full glory in the Big Table.%
\begin{equation*}
\end{equation*}

\textbf{10. Reflexive polytopes.}

M.Kreuzer and H.Skarke have classified reflexive polyhedra in four
dimensions . They found 473 800 652 of them. \ For each of them is
associated a Laurent polynomial \ $S$ \ in four variables \ Then a solution
to a Calabi-Yau differential equation is constructed with coefficients \ $%
A_{n}=$constant term$(S^{n}).$

We will show the idea with an example in dimension two where there are only
16 reflexive polytopes

\textbf{Example. }Consider the polygon with four vertices%
\begin{equation*}
\left( 
\begin{tabular}{llll}
2 & 0 & -1 & -1 \\ 
-1 & 1 & 0 & -1%
\end{tabular}%
\right)
\end{equation*}%
with associated Laurent polynom%
\begin{equation*}
S=\frac{x^{2}}{y}+y+\frac{1}{x}+\frac{1}{xy}
\end{equation*}%
Then 
\begin{equation*}
A_{n}=\text{c.t.}(S^{n})=\sum_{k}\binom{n}{3k}\binom{3k}{k}\binom{2k}{n-4k}
\end{equation*}%
Using "Zeilberger" we get a differential equation of order six which factors
into a huge left factor of order four and a right factor

\begin{equation*}
5\theta ^{2}+x\theta (11\theta -1)+6x^{2}\theta ^{2}+x^{3}\theta (13\theta
-9)-x^{4}(298\theta ^{2}+1636\theta +960)
\end{equation*}%
\begin{equation*}
-4x^{5}(726\theta ^{2}+3510\theta +2391)-8x^{6}(917\theta ^{2}+4752\theta
+3519)
\end{equation*}%
\begin{equation*}
-32x^{7}(256\theta ^{2}+1372\theta +1077)-32x^{8}(327\theta ^{2}+810\theta
+528)
\end{equation*}%
\begin{equation*}
-576x^{9}(54\theta ^{2}+78\theta +29)-128x^{10}(\theta +1)(355\theta +417)
\end{equation*}%
\begin{equation*}
-23808x^{11}(\theta +1)(\theta +2)
\end{equation*}%
with solution%
\begin{equation*}
y_{0}=1+12x^{4}+60x^{5}+420x^{8}+...
\end{equation*}%
This example shows that a very simple polytope can give a differential
equation of very high degree.

In dimension four we consider

\textbf{Example \#325. }Consider the polytope with 13 vertices giving the
Laurent polynomial%
\begin{equation*}
S=\frac{1}{x}+\frac{y}{x}+\frac{x}{y}+\frac{z}{x}+\frac{x}{z}+\frac{yz}{x}+%
\frac{x}{yz}+\frac{t}{x}(1+y+yz)+\frac{x}{t}(1+\frac{1}{y}+\frac{1}{yz})
\end{equation*}%
Then there are only even terms so we compute%
\begin{equation*}
A_{n}=\text{c.t.}(S^{2n})
\end{equation*}%
A direct approach by expanding the powers of \ $S$ \ and then take the
constant term costs a lot of computer time. Here we will eliminate \ $t$ \
and thus reducing the computer time by a factor 500-1000. Let%
\begin{equation*}
u:=\frac{1}{x}+\frac{y}{x}+\frac{x}{y}+\frac{z}{x}+\frac{x}{z}+\frac{yz}{x}+%
\frac{x}{yz}
\end{equation*}%
\begin{equation*}
p=\frac{t}{x}(1+y+yz)
\end{equation*}%
\begin{equation*}
q=\frac{x}{t}(1+\frac{1}{y}+\frac{1}{yz})
\end{equation*}%
\begin{equation*}
v=pq=(1+y+yz)(1+\frac{1}{y}+\frac{1}{yz})
\end{equation*}%
Then we have 
\begin{equation*}
A_{n}=\text{c.t.(}\sum_{i+j+k=2n}\frac{(2n)!}{i!j!k!}u^{i}p^{j}q^{k})
\end{equation*}%
To get rid of \ $t$ \ we need \ $j=k.$ It results%
\begin{equation*}
A_{n}=\text{c.t.}(\sum_{i+2j=2n}\frac{(2n)!}{i!j!^{2}}u^{i}(pq)^{j})=\text{%
c.t.(}\sum_{n=0}^{n}\frac{(2n)!}{j!^{2}(2n-2j)!}u^{2n-2j}v^{j})
\end{equation*}%
To find the equation for 
\begin{equation*}
y_{0}=\sum_{n=0}^{\theta }A_{n}x^{n}
\end{equation*}%
we need 30 coefficients which are computed in about five minutes on a
laptop. \ Arne Meurman has computed an explicit formula for \ $A_{n}$ \
summing over eight indices with complicated summation limits, making it not
very \ useful.The equation is%
\begin{equation*}
19^{2}\theta ^{4}-19x(4333\theta ^{4}+6212\theta ^{3}+4778\theta
^{2}+1672\theta +228)
\end{equation*}%
\begin{equation*}
+x^{2}(4307495\theta ^{4}+7600484\theta ^{3}+6216406\theta
^{2}+2802424\theta +530556)
\end{equation*}%
\begin{equation*}
-x^{3}(93739369\theta ^{4}+213316800\theta ^{3}+236037196\theta
^{2}+125748612\theta +25260804)
\end{equation*}%
\begin{equation*}
+x^{4}(240813800\theta ^{4}+778529200\theta ^{3}+1041447759\theta
^{2}+631802809\theta +138510993)
\end{equation*}%
\begin{equation*}
-2^{2}\cdot 409x^{5}(\theta +1)(2851324\theta ^{3}+100355\theta
^{2}+11221241\theta +3481470)
\end{equation*}%
\begin{equation*}
+2^{2}\cdot 3^{2}\cdot 19^{2}\cdot 409^{2}x^{6}(\theta +1)(\theta
+2)(2\theta +1)(2\theta +5)
\end{equation*}%
\begin{equation*}
\end{equation*}

Let us consider another example,with 18 vertices, v18.16766 in the notation
of Batyrev-Kreuzer [12] with%
\begin{equation*}
S=x(1+t+yt+zt+yzt^{2})+\frac{1}{x}(1+\frac{1}{t}+\frac{1}{yzt^{2}})
\end{equation*}%
\begin{equation*}
+y+\frac{1}{y}+z+\frac{1}{z}+\frac{1}{t}+yt+\frac{1}{yt}+zt+\frac{1}{zt}+yzt
\end{equation*}%
Kreuzer computed 60 coefficients with brute force (time 760 hours). Using
the the method above eliminating \ $x$ \ it took only 71 minutes. Finally
Duco van Straten, Pavel Metelitsyn and Elmar Sch\"{o}mer using modular
arithmetic computed 272 coefficients. They succeded in finding a
differential equation of order \ $6$ \ and degree \ $25.$It is not MUM and
does not factor. This example leaves some doubts about the reflexive
polytopes.

So far we know polytopes giving the coefficients of the following equations
in the Big Table (see

1-14, 16, 24, 25, 26, 29, 42, 51, 70, $\sim $101, 185, 188, 206, 209, 214,
218, 287, 308, 309, 324,325, 326, 327

Only the last four are new. 
\begin{equation*}
\end{equation*}

\textbf{11. Other equations.}

\textbf{11.1. Bessel moments.}

Consider the Bessel moments (see D.Bailey, J.Borwein, D.Broadhurst and
M.L.Glasser [9]) 
\begin{equation*}
c_{m,k}=\dint\limits_{0}^{\infty }x^{k}K_{0}(x)^{m}dx
\end{equation*}%
Here \ $K_{0}(x)$ \ is a certain Bessel function that conveniently can be
defined by%
\begin{equation*}
K_{0}(x)=\dint\limits_{0}^{\infty }e^{-x\cosh (t)}dt
\end{equation*}%
This leads to another representation (in Ising theory)%
\begin{equation*}
c_{m,k}=\frac{k!}{2^{m}}\dint\limits_{0}^{\infty
}...\dint\limits_{0}^{\infty }\frac{dx_{1}...dx_{m}}{(\cosh
(x_{1})+...+\cosh (x_{m}))^{k+1}}
\end{equation*}%
Let 
\begin{equation*}
d_{n}=\frac{15^{2n}}{n!^{2}}c_{5,2n-1}
\end{equation*}%
and%
\begin{equation*}
Y_{0}=\sum_{n=0}^{\infty }d_{n}x^{n}
\end{equation*}%
Then $\ Y_{0}$ \ satisfies the differential equation%
\begin{equation*}
\theta ^{2}(\theta -1)^{2}-4x\theta ^{2}(259\theta
^{2}+26)+3600x^{2}(35\theta ^{4}+70\theta ^{3}+63\theta ^{2}+28\theta
+5)-3240000x^{3}(\theta +1)^{4}
\end{equation*}%
The shape of the last term suggests that converting to \ $x=\infty $ \ could
give a Calabi-Yau equation. Indeed \ $\theta \longrightarrow -\theta -1$ \
and \ $x\longrightarrow 900x^{-1}$ \ gives the equation

\begin{equation*}
\theta ^{4}-x(35\theta ^{4}+70\theta ^{3}+63\theta ^{2}+28\theta +5)
\end{equation*}%
\begin{equation*}
+x^{2}(\theta +1)^{2}(259\theta ^{2}+518\theta +285)-225x^{3}(\theta
+1)^{2}(\theta +2)^{2},
\end{equation*}%
which we recognize as \# 34 (found by H.Verrill [19]) with solution \ $%
y_{0}=\dsum\limits_{n=0}^{\theta }A_{n}x^{n}$ \ with%
\begin{equation*}
A_{n}=\sum_{i+j+k+l+m=n}\left( \frac{n!}{i!j!k!l!m!}\right) ^{2}
\end{equation*}%
Similarly%
\begin{equation*}
d_{n}=\frac{48^{2n}}{n!^{2}}c_{6,2n-1}
\end{equation*}%
leads to a differential equation whose mirror at infinity is the 5-th order
equation \#130 (also found by Verrill) with coefficients%
\begin{equation*}
A_{n}=\sum_{i+j+k+l+m+s=n}\left( \frac{n!}{i!j!k!l!m!s!}\right) ^{2}
\end{equation*}

\textbf{11.2. Differential equations coming from combinatorics.}

There are probably many Calabi-Yau equations whose coefficients count
something in combinatorics. The last equation \# 366 in the Big Table comes
from counting random walks in \ $\mathbb{Z}^{4}$ (communicated to me by
T.Guttmann, [17] )

\begin{equation*}
\theta ^{4}+x\theta (39\theta ^{3}-30\theta ^{2}-19\theta -4)
\end{equation*}%
\begin{equation*}
+2x^{2}(16\theta ^{4}-1070\theta ^{3}-1057\theta ^{2}-676\theta -192)
\end{equation*}%
\begin{equation*}
-2^{2}3^{2}x^{3}(3\theta +2)(171\theta ^{3}+566\theta ^{2}+600\theta +316)
\end{equation*}%
\begin{equation*}
-2^{5}3^{3}x^{4}(384\theta ^{4}+1542\theta ^{3}+2635\theta ^{2}+2173\theta
+702)
\end{equation*}%
\begin{equation*}
-2^{6}3^{3}x^{5}(\theta +1)(1393\theta ^{3}+5571\theta ^{2}+8378\theta +4584)
\end{equation*}%
\begin{equation*}
-2^{10}3^{5}x^{6}(\theta +1)(\theta +2)(31\theta ^{2}+105\theta +98)
\end{equation*}%
\begin{equation*}
-2^{12}3^{7}x^{7}(\theta +1)(\theta +2)^{2}(\theta +3)
\end{equation*}%
It was found by computing 40 coefficients. An explicit formula for the
coefficients is unknown.

Also \# 16 comes from combinatorics, see [16].%
\begin{equation*}
\end{equation*}

\textbf{12. Some remarks.}

Let \ $p$ \ be a prime. Expand in base \ $p$%
\begin{equation*}
n=n_{0}+n_{1}p+n_{2}p^{2}+...
\end{equation*}%
\begin{equation*}
k=k_{0}+k_{1}p+k_{2}p^{2}+...
\end{equation*}%
Then it is wellknown that%
\begin{equation*}
\binom{n}{k}\equiv \binom{n_{0}}{k_{0}}\binom{n_{1}}{k_{1}}\binom{n_{2}}{%
k_{2}}....\func{mod}\text{ }p
\end{equation*}%
Kira Samol and Duco van Straten have found that a similar congruences (Dwork
congruences) are valid for the coefficients \ $A_{n}$ of most Calabi-Yau
equations, namely%
\begin{equation*}
A_{n}\equiv A_{n_{0}}A_{n_{1}}A_{n_{2}}...\func{mod}\text{ }p
\end{equation*}%
It seems to be valid also for very complicated coefficients, like \#264, 274
and also for \# 366 where we have no formula for \ $A_{n}.$ This property is
not preserved under equivalence transformations so it is not surprising that
it is not valid for some pullbacks of fifth order equations.

Hundreds of binomial identities resulted from the search of Calabi-Yau
differential equations. We give only a small sample of the simplest
identities%
\begin{equation*}
\sum_{k}\binom{n}{k}^{2}\binom{3k}{2n}=\sum_{k}\binom{n}{k}^{2}\binom{2k}{k}
\end{equation*}%
\begin{equation*}
\sum_{k}(-1)^{n+k}\binom{n}{k}\binom{n+k}{n}^{2}=\sum_{k}\binom{n}{k}^{2}%
\binom{n+k}{n}
\end{equation*}%
\begin{equation*}
\sum_{k=0}^{n}(-1)^{k+l}\binom{n}{k}\binom{n}{2l-k}=\binom{n}{l}
\end{equation*}%
\begin{equation*}
\sum_{l=0}^{n}(-1)^{n+k+l}\binom{n}{l}\binom{2l}{n-k}\binom{2n-2l}{k}=2^{n}%
\binom{n}{k}
\end{equation*}%
\begin{equation*}
\end{equation*}

\textbf{Acknowledgements.}

First of all I want to thank my collaborators Christian Krattenthaler, Duco
van Straten and Wadim Zudilin for working with me during many years. Further
I thank Christian van Enckevort, Jesus Guillera, Tony Guttmann, Max Kreuzer,
Arne Meurman, Peter Paule, Carsten Schneider, Helena Verrill and Don Zagier
for various contributions. Finally I thank Maple for providing me with
numerous recursion formulas when "Zeilberger" is applied to sums that are
identical zero.%
\begin{equation*}
\end{equation*}

\textbf{References.}

\textbf{1. }G.Almkvist, Str\"{a}ngar i m\aa nsken (in Swedish),I, Normat 
\textbf{51 }(2003), no.1, 22-33, II, Normat \textbf{51 }(2003), no.2,63-79.

\textbf{2.}G.Almkvist, W.Zudilin, Differential equations, mirror maps and
zeta values, in: Mirror Symmetry V, N.Yui,

\ \ S.-T. Yau, and J.D. Lewis (eds), AMS/IP Stud. Adv. Math. 38
(International Press \& Amer. Math. Soc.,

\ \ Providence,RI, 2007) 481-515; arXiv: math/0402386

\textbf{3.}G.Almkvist, C.van Enckevort, D.van Straten, W.Zudilin, Tables of
Calabi-Yau equations, arXiv:

\ \ math/0507430.

\textbf{4.}G.Almkvist, D.van Straten, W.Zudilin, Generalizations of
Clausen's formula and algebraic trnsformations

\ \ of Calabi-Yau differential equations, in preparation.

\textbf{5. }G.Almkvist, C.Krattenthaler, Some harmonic sums related to
Calabi-Yau differential equations, in preparation.

\textbf{6.} G.Almkvist, Calabi-Yau differential equations of degree 2 and 3
and Yifan Yang's pullback, AG/0612215.

\textbf{7. }\ G.Almkvist, 5-th order differential equations related to
Calabi-Yau differential equations, AG/0703261

\textbf{8. }G.Almkvist, Ap\'{e}ry, Bessel, Calabi-Yau and Verrill,
CA/08081480

\textbf{9.}D.H.Bailey. J.M.Borwein, D.Broadhurst, M.L.Glasser, Elliptic
integral evaluations of Bessel moments, arXiv:

\ \ hep-th/0801089.

\textbf{10. }V.V.Batyrev, D. van Straten, Generalized hypergeometric
functions and rational curves on Calabi-Yau

\ \ \ complete intersections in toric varieties, Commun. Math. Phys. \textbf{%
168} (1995), 493-533.

\textbf{11. }V.V.Batyrev, I.Ciocan-Fontanine, B.Kim,D. van Straten, Conifold
transitions and mirror symmetry for

\ \ Calabi-Yau complete intersections in Grassmannians, Nuclear Phys. B 
\textbf{514 }(1998), no.3, 640-666, AG/9710022

\textbf{12. }V.V.Batyrev, M.Kreuzer, Constructing new Calabi-Yau 3-folds and
their mirror conifold transitions, AG/08023376.

\textbf{13.}J.M.Borwein, B.Salvy, A proof of a recursion for Bessel moments,
inria-00152799

\textbf{14. }C.Doran, J.Morgan, Mirror symmetry and integral variations of
Hodge structure underlying one parameter

\ \ \ \ \ families of Calabi-Yau threefolds, AG/0505272.

\textbf{15. }J.Guillera, About a new kind of Ramanujan-type series,
Experimental Math. \textbf{12} (2003), no.4, 507-510.

\textbf{16. }M.L.Glasser, A.J.Guttmann, Lattice Green functions (at 0) for
the 4D hypercubic lattice, J. Phys. A Math. Gen. vol \textbf{27} (1994),
7011-14.

\textbf{17. }T.Guttmann, Lattice Green functions in 4 and 5 dimensions and
Calabi-Yau equations, in preparation.

\textbf{18. }P.Paule, C.Schneider, Complete proofs of a new family of
harmonic number identities, Preprint, RISC, 2002

\textbf{19. }E.A.R\"{o}dland, The Pfaffian Calabi-Yau, its mirror, and their
link to Grassmannian G(2,7), Composito Math. \textbf{122 }(2000), no.2,
135-149, AG/9801092

\textbf{20. }H.Verrill,Sums of squares of binomial coefficients, with
applications to Picard-Fuchs differential

\ \ equations, math.CO/0407327.

\textbf{21. }Y.Yang, personal communication.

\textbf{22. }D.Zagier, Integral solutions of Ap\'{e}ry-like recurrence
equations,

\textbf{22.}W.Zudilin, Binomial sums related to rational approximations of $%
\varsigma (4),$Mat. Zametki 75 (2004), 637-640,

\ \ \ \ English transl. Math. Notes \textbf{75} (2004), 594-597,
math.Ca/031196.

\ 

\end{document}